\documentclass[11pt]{article}
\usepackage{amsfonts,amsmath,amssymb,amsthm, amscd}
\usepackage[dvips]{epsfig,graphics}

\numberwithin{equation}{section}
\newtheorem{theorem}{Theorem}[section]

\newtheorem{cor}[theorem]{Corollary}
\newtheorem{conj}[theorem]{Conjecture}

\theoremstyle{definition}
\newtheorem{definition}[theorem]{Definition}
\newtheorem{example}[theorem]{Example}
\newtheorem{remark}[theorem]{Remark}

\def\<{{\langle}}
\def\>{{\rangle}}

\def\Z{\mathbb Z}
\def\Q{\mathbb Q}

\def\S{{\mathbb S}}

\def\t{\tau}
\def\w{\omega}

\def\ni{\noindent} 
\begin{document}

\title{Knot Groups with Many Killers}

\author{Daniel S. Silver \thanks{Partially supported by NSF grant
DMS-0706798.}\and Wilbur Whitten\and{Susan G. Williams \thanks{Partially supported by NSF grant DMS-0706798}} }

\maketitle %{\setlength{\linewidth}{2in}

\begin{abstract}   \noindent The group of any nontrivial torus knot, hyperbolic 2-bridge knot, or hyperbolic knot with unknotting number one contains infinitely many elements, none the automorphic image of another, such that each normally generates the group.  \end{abstract}

 \begin{footnote}{Mathematics Subject Classification:  
Primary 57M25; secondary 20E34.}\end{footnote}

%%%%%%%%%%%%%%%%%%%%%%%%%%%%%% 1. INTRODUCTION %%%%%%%%%%%%%
\section{Introduction} \label{Section 1} Let $k$ be a knot in $\S^3$ with group $\pi_k=\pi_1(\S^3 \setminus k, *)$. For convenience, we choose the basepoint $*$ on the boundary of a regular neighborhood $k \times {\mathbb D}^2$.  A {\it meridian} is an element $\mu \in \pi_k$ represented by a curve freely homotopic to $* \times \partial {\mathbb D}^2$. We orient $k$ and insist that the curve have linking number $+1$ with the knot. Then $\mu$ is well defined up to conjugation.  

Following \cite{simon}, we say that an element of a group is a {\it killer} if the group modulo the element is trivial; in other words, the element normally generates the group. Obviously the image of a killer under any automorphism of the group is a killer. 
It is well known that the meridian of a knot group is a killer. In \cite{tsau1}, C.M. Tsau gave an example of a knot group killer that is not the automorphic image of the meridian. He called such an element a ``nonalgebraic killer." We prefer a less violent term,  ``pseudo-meridian." 

\begin{definition} A {\it pseudo-meridian} of a knot group is an element that normally generates the group but is not an automorphic image of the meridian. Two pseudo-meridians are {\it equivalent} if one is the automorphic image of the other.  \end{definition}

\begin{theorem}\label{mainthm} Let $k$ be a nontrivial $2$-bridge knot or  torus knot, or  a hyperbolic knot with unknotting number one. Then $\pi_k$ contains infinitely many pairwise nonconjugate pseudo-meridians. \end{theorem}

It is well known that the group of automorphisms of $\pi_k$  modulo inner automorphisms is finite. This is a consequence of Mostow Rigidity when $k$ is hyperbolic (see \cite{thurston}, Chapter 5). When $k$ is a torus knot, the group has order 2 (see \cite{kaw}).
Hence we obtain the following.

\begin{cor} \label{nonequivalent}  Let $k$ be a nontrivial $2$-bridge knot or  torus knot, or a  hyperbolic knot with unknotting number one. Then $\pi_k$ contains infinitely many pairwise nonequivalent pseudo-meridians. \end{cor}

Our interest in pseudo-meridians is motivated largely by a conjecture of Jonathan Simon \cite{kirby}:   a sequence $\pi_{k_1} \to \pi_{k_2} \to \ldots \to \pi_{k_n}\to \ldots $ of knot group epimorphisms must be isomorphisms, for sufficiently large $n$. Much has been written about Simon's conjecture (see \cite{bbrw} for numerous references). It was recently proven for all 2-bridge knots \cite{bbrw}.

Given a knot group epimorphism, the image of a meridian is either a meridian or a pseudo-meridian. Examples for which the image is a pseudo-meridian are known \cite{jl}, \cite{sw}.  It is not known whether an epimorphism mapping a meridian to meridian must exist.

\section{Proof of Theorem \ref{mainthm}.}  {\sl Assume that $k$ is a nontrivial 2-bridge knot.}  From a 2-bridge diagram for $k$, we obtain a group presentation  for $\pi_k$ having the form $\< x, y \mid r\>$, where $x$ and $y$ are meridians corresponding to arcs of a 2-bridge diagram for $k$. Without loss of generality, we can assume that the exponent sum of $y$ in $r$ is $+1$ (and consequently, that of $x$ is $-1$).  We introduce a new generator $a$ and defining relation 
$y = ax$, and use the relation to eliminate $y$ from the presentation. 
We obtain a presentation of the form
$\pi_k = \<x, a \mid \prod x^{-k_i}a^{\epsilon_i}x^{k_i}\>,$
where $\sum \epsilon_i =0$. 

Consider the elements $\mu_n = x(yx^{-1})^n= xa^n,$ where $n \ge 0$. 
Killing $\mu_n$ introduces the relation $x = a^{-n}$. The relator above becomes $\prod a^{-nk_i} a^{\epsilon_i} a^{nk_i}$, which is simply $a$.
Since $x = a^{-n}$, killing $\mu_n$ also kills $x$. Hence $\mu_n$ is 
a killer, for any $n \ge 0$. We will show that the elements $\mu_n$ are pairwise nonconjugate, for sufficiently large $n$. 

The group $\pi$ admits a nonabelian parabolic representation $\rho: \pi_k \to {\rm SL}_2{\mathbb C}$: 
\begin{equation} \label{matrices} x \mapsto X = \begin{pmatrix} 1& 1 \\ 0 & 1 \end{pmatrix},\quad y \mapsto Y = \begin{pmatrix}1& 0 \\ \w & 1 \end{pmatrix},\end{equation}
where $\w$ is any root of a nontrivial polynomial $\Phi(w)$ defined in \cite{riley}. 
Denote the trace of $\rho(\mu_n) = X(YX^{-1})^n$ by $\t_n$. 
Some values of $\t_n$ are given below. 

$$\begin{matrix} n & \t_n \\
-- & -- \\
0 & 2 \\
1 & 2 \\
2 & 2(1-\w) \\
3 & 2(1-3\w+\w^2)\\
4 & 2(-1+6\w -5\w^2+\w^3) \end{matrix}$$
(The values $\t_0$ and $\t_1$ are both 2 for an obvious reason: $\mu_0 = x$ while $\mu_1 = xyx^{-1}$ is conjugate to $x$.)

One checks that $\t_n$ satisfies a linear recurrence relation:
$$\t_n = (2-\w) \t_{n-1} - \t_{n-2},\quad \t_0 = \t_1 =2.$$
The characteristic equation of the recurrence is $x^2 + (\w-2) x +1 =0$.
The roots of the equation are $\l, \l^{-1}$ such that $\l + \l^{-1} = 2-\w$. 
It is a straightforward matter to solve: 
$$\t_n = 2 \cdot {  {\l^{n -\frac12} + \l^{\frac12-n} } \over {\l^{\frac12}+ \l^{-\frac12} }}.$$
If $|\l|\ne 1$, then the absolute value of $\t_n$ grows exponentially, and hence the elements $\mu_n$ are pairwise nonconjugate, for $n$ sufficiently large. The claim is established in this case.

Suppose that $\l = e^{i\theta}$, for some $\theta$. The 
expression for $\t_n$ simplifies:
$$\t_n = 2 \cdot {  \cos((n-\frac12)\theta) \over \cos(\frac\theta 2)}.$$
Either $\t_n$ is non-repeating, and the claim is proved, or else $\theta$ is rationally related to $\pi$. 

Suppose that $\theta$ is rationally related to $\pi$; that is, suppose that $\l = e^{i\theta}$, where $\theta \in \Q\pi$. 
We can use the relation $\l + \l^{-1} = 2-\w$ to solve for $\w$: 
$$\w = 2(1-\cos \theta) = 4 \sin^2 (\tfrac\theta 2)>0.$$

 If $k$ is a torus knot, then the coefficients of $\Phi(w)$ are all positive (an easy consequence of the form of the 2-bridge knot group presentation). Since not all the roots of such a polynomial can be positive, we can choose a root $\w$ away from the positive real axis. Then the traces $\t_n, \ n \ge 1,$ are non-repeating. (In fact, all roots of $\Phi(w)$ are real and negative. We don't need this here.)
 
If $k$ is hyperbolic, then for some choice of $\w$, $\rho$ is the representation corresponding to the hyperbolic structure. This $\w$ cannot be real since $\pi_k$ is a Kleinian group of finite covolume (see Exercise 1.3, No. 1 of \cite{reid}, for example). Again the traces $\t_n, \ n \ge 1,$ are non-repeating. A 2-bridge knot cannot be a satellite, so Theorem \ref{mainthm} is proved for 2-bridge knots. 

{\sl Next suppose that $k$ is a hyperbolic knot with unknotting number one.} There exists a crossing in some diagram of $k$ such that changing the crossing results in a diagram of the trivial knot. We regard this crossing as a 4-valent vertex with meridianal generators $x, y, z, w$, as in Figure 1. Suppose that $k$ has a left-hand crossing here. (The right-hand case is similar.) Then $\pi_k$ has a presentation with relations $xy=zw$ and $x=w$, together with Wirtinger generators and relations for the other arcs and crossings.  Replacing $x=w$ with $y=z$ gives a presentation of $\Z$, the group of the trivial knot.

\begin{figure}
\begin{center}
\includegraphics[height=1 in]{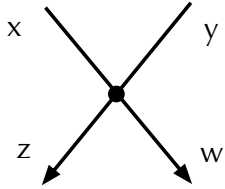}
\caption{Crossing and generators}
\label{crossing}
\end{center}
\end{figure}

Let $\mu_n = x a^{-n}$, where $a = y x^{-1}$ and $n$ is an arbitrary integer. A presentation for the quotient of $\pi_k$ obtained by killing $\mu_n$ can be obtained as follows: set $x$ and $y$ equal to $a^n$ and $a^{n+1}$, respectively, leaving other Wirtinger generators unchanged. The relations at the crossing in Figure 1 force $z=a^{n+1}$ and $w=a^n$.   Applying the same procedure in the given presentation of the unknot has the same result, so the two quotient groups are equal.  But the quotient of $\Z$ obtained by setting one generator equal to $a^n$ and another equal to $a^{n+1}$ is trivial.   Hence $\mu_n$ is a killer, for any $n$. 

The elements $x, y$ do not commute in $\pi_k$. (Otherwise, $y=z$ in $\pi_k$. As above, changing the crossing would not affect the group.
However the group would be trivial. )  By \cite{riley}, there exists a parabolic representation $\rho: \pi_k \to {\rm SL}_2 {\mathbb C}$ such that $\rho(x), \rho(y)$ have the form in equation (\ref{matrices}). 
The same argument as in the case of 2-bridge knots shows that the $\mu_n$ are pairwise nonconjugate.

{\it Finally, assume that $k$ is a nontrivial $(p, q)$-torus knot.} Its group has a presentation of the form

$$ \pi_k =\<u, v \mid u^p=v^q\>. $$
We can assume without any loss of generality that $p, q$ are relatively prime integers with $p > q>0$. 

There is an epimorphism $\chi: \pi_k \to \Z$ with $\chi(u) = q$ and $\chi(v) = p$. We can find integers $r$, $s$ with 
\begin{equation}\label{rs} rp+sq=1, |r|<q, |s|<p \end{equation}
Then $\chi(u^sv^r) = 1$. Notice that by replacing $r$ and $s$ with $r-q$ and $s+p$, if necessary, we can assume
that $s$ is positive in condition \ref{rs}. 
Clearly, $s=1$ is impossible. Hence $s>1$.  

Introduce new generators $x, a$ and defining relations $x=u^s v^r, a=ux^{-q}$.  Then

$$\pi_k= \<x, u, v, a \mid u^p=v^q, x=u^sv^r, u=ax^q\>$$

Using the relations above, eliminate $u$. We have

\begin{equation}\label{**}\pi_k= \<x, v, a \mid (ax^q)^p= v^q, (ax^q)^{-s}x = v^r\> \end{equation}

We claim that the element $\mu_n=xa^{-n}$ is a killer, for any integer $n$. Setting $x=a^n$ in the relations in presentation (\ref{**}) gives
$$a^{p(nq+1)}=v^q, a^{-s(nq+1)+n}=v^r,$$
and so 
$$v^{rq}= a^{rp(nq+1)}=a^{-sq(nq+1)+nq}.$$
Since $rp+sq=1$, we have $a^{nq+1}=a^{nq}$. Thus $a=1$, and so $x=1$. It follows immediately that $v=1$. 

To show that the $\mu_n$ are pairwise nonconjugate, it suffices to show their images are pairwise nonconjugate in the quotient
$$\<u, v \mid u^p=v^q=1\> = \Z/(p)*\Z/(q).$$
Rewriting $xa^{-n}$ in terms of $u$ and $v$ yields 
$$\mu_n=u^sv^r[u(u^sv^r)^{-q}]^{-n}=u^sv^r[(u^sv^r)^{q}u^{-1}]^{n}.$$

Since $s>1$,  cyclic reduction of this word won't reduce the number of occurrences of $v^{\pm1}$, which is a monotone function of $n$. Hence the killers $\mu_n$ are pairwise nonconjugate.

%note that the relations in presentation (\ref{**}) imply 

%
%$$ (ax^q)^{rp}=v^{rq}, [(ax^q)^{-s}x]^q= v^{rq}.$$
%Hence $(ax^q)^{rp} = [(ax^q)^{-s}x]^q.$ Letting $x=a^{-n}$, we have
%$(a^{-qr+1})^{sq}(a^{-qr+1})^{rp} = a^{-nq}$.  Since $rp+ sq=1$, we have 
%$a^{-qn+1}=a^{-qn}$. Thus $a=1$, and so $x=1$. It follows immediately that $v=1$. 

%
%In view of Lemma \ref{nonequivalent},  it suffices to show the elements are pairwise nonconjugate in the quotient
%$$\<u, v \mid u^p=v^q=1\> = \Z/(p)*\Z/(q).$$
%Rewriting $xa^n$  in terms of $u, v$ yields:

%$$ u^sv^r[u(v^{-r} u^{-s})^q]^n. $$
%After cyclic reduction, we have 
%$$[u(v^{-r} u^{-s})^q]^{n-1}u(v^{-r}u^{-s})^{q-1}.$$
%Since $s>1$, further cyclic reduction won't reduce the number of occurrences of $v$, which is a monotone function of $n$. Hence the killers $\mu_n$ are pairwise nonconjugate. 

The proof of Theorem \ref{mainthm} is complete. 

\begin{remark} Consider a knot diagram with a crossing, described by Figure 1, with the property that the elements $\mu_n = x(yx^{-1})^n$ are pairwise nonequivalent killers for arbitrarily large $n$. Introducing any number of full twists in the two arcs corresponding to $z, w$ preserves this property, provided that the resulting knot is hyperbolic. For letting $x=a^n$ and $y=a^{n+1}$, as in the proof of Theorem \ref{mainthm}, results in the same trivial quotient group. The above argument for hyperbolic knots applies. Hence twisting in this manner produces many new examples of knot groups with infinitely many pairwise nonequivalent killers. \end{remark} 

\section{Examples and conjecture.}

%We conclude with two examples (\ref{ex1}, \ref{ex2}) and a conjecture (\ref{conj}). 

\begin{example} \label{ex1} Consider the diagram for the hyperbolic knot $k= 8_{20}$ in Figure 2, with certain Wirtinger generators indicated. Changing the crossing involving $x$ and $y$ produces a diagram of the unknot. Hence by the proof of Theorem \ref{mainthm}, the elements $\mu_n = x(yx^{-1})^n$ are nonequivalent pseudo-meridians, for sufficiently large $n$ (in fact, for all $n$). 

\begin{figure}
\begin{center}
\includegraphics[height=2 in]{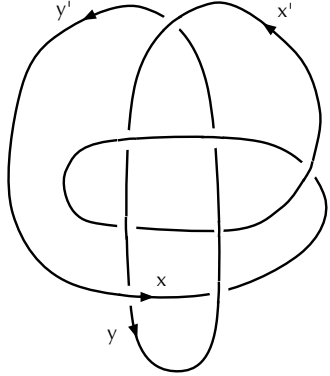}
\caption{The knot $8_{20}$}
\label{knot}
\end{center}
\end{figure}

From the proof of Theorem \ref{mainthm}, one might suspect that 
the elements $x(yx^{-1})^n$ are killers of any knot group whenever $x$ and $y$ are noncommuting meridians. (If so, then the conclusion of Theorem \ref{mainthm} would follow for any hyperbolic knot.) However, this is not the case for this example. If we choose meridians
$x'$ and $y'$, then this word, according to GAP \cite{gap},  is not a killer for $n=2$. 

\end{example}

\begin{example} \label{ex2} Let $\tilde k$ be a hyperbolic knot embedded in a standard solid torus $V$ in such a way that $V\setminus \tilde k$ is a hyperbolic manifold. The pair $(V, \tilde k)$ is  represented by a link diagram, one component representing $\tilde k$ and the other a meridian of $V$. Assume that changing some crossing unknots $\tilde k$ in $V$ (equivalently, it produces a trivial 2-component link diagram). Twist knots provide simple examples. 

Let $\hat k$ be any knot. Form the satellite knot $k$ with companion $\hat k$ and pattern $(V, \tilde k)$ (see \cite{bz}, for example). Changing a crossing of $\tilde k$ that unknots it in $V$ also unknots the satellite knot $k$. As in the proof of Theorem \ref{mainthm}, we produce a sequence $\mu_n = x (yx^{-1})^n$ of elements that are killers for $\pi_{\tilde k}$. They are also killers for $\pi_k$. This can be seen from the link diagram for the pattern:  As in the proof of Theorem \ref{mainthm}, the quotient group of $\pi_1(V \setminus \tilde k)$ we get by setting $x=a^n$ and $y=a^{n+1}$ is unaffected if first we change the crossing and then make the substitution. But changing the crossing produces a trivial 2-component link, and then the substitution kills $x$. Since $\pi_1(V\setminus \tilde k)$ is a subgroup of $\pi_k$ and since $x$ normally generates $\pi_k$, the elements $\mu_n$ are killers of $\pi_k$. 

If $\tilde k$ is hyperbolic, then the argument in the proof of Theorem \ref{mainthm} shows that, for $n$ sufficiently large, no $\mu_n$ is the image of another under an automorphism of $\pi_{\tilde k}$. 

Since $V \setminus \tilde k$ is a hyperbolic manifold, and since the characteristic submanifold of the exterior of $k$ is unique up to ambient isotopy   \cite{jaco}, it follows that any autormorphism of $\pi_k$ can be realized by a homeomorphism that leaves $\partial V$ invariant, mapping the longitude $\lambda \in \pi_1(\partial V)$ to itself. The restriction to $V\setminus \tilde k$ induces an automorphism of  $\pi_1(V\setminus \tilde k)$ which induces an automorphism of $\pi_{\tilde k}\cong \pi_1(V \setminus \tilde k)/\<\<\lambda\>\>$. Hence the pseudo-meridians $\mu_n$ for $\pi_{\tilde k}$ that are pairwise nonequivalent are also 
pairwise nonequivalent pseudo-meridians for $\pi_k$. 

\end{example} 

\begin{conj}\label{conj} Every nontrivial knot group has infinitely many nonequivalent pseudo-meridians. 

\end{conj}
%%%%%%%%%%%%%%%%%%%%%%%%%

 \bigskip

\ni Department of Mathematics and Statistics,\\
\ni University of South Alabama\\ Mobile, AL 36688 USA\\
\ni Email: silver@jaguar1.usouthal.edu\\ 

\ni 1620 Cottontown Road,\\
\ni Forest, VA 24551 USA\\
\ni Email: bjwcw@aol.com \\ 

\ni Department of Mathematics and Statistics,\\
\ni University of South Alabama\\ Mobile, AL 36688 USA\\
\ni Email: swilliam@jaguar1.usouthal.edu
\end{document}